\documentclass[11pt]{amsart}
\usepackage[T1]{fontenc}
\usepackage{amssymb}
\usepackage{amsmath}
\usepackage{fancyhdr}
\usepackage[british]{babel}
\usepackage{geometry,mathtools}
\usepackage{enumitem}
\usepackage{algpseudocode}
\usepackage{dsfont}
\usepackage{centernot}
\usepackage{xstring}
\usepackage{colortbl}
\usepackage[section]{algorithm}
\usepackage{graphicx}
\usepackage[font={footnotesize}]{caption}
\usepackage[usenames,dvipsnames,table]{xcolor}
\usepackage{pstricks,tikz}
\usetikzlibrary{patterns}
\usepackage[h]{esvect}
\usepackage[
bookmarksopen=true,
bookmarksopenlevel=1,
colorlinks=true,
linkcolor=darkblue,
linktoc=page,
citecolor=darkblue,
]{hyperref}
\usepackage{bbm}
\definecolor{darkblue}{rgb}{0,0,0.5}

\newtheorem{theorem}[algorithm]{Theorem}

\newtheorem{cor}[algorithm]{Corollary}
\newtheorem{observe}[algorithm]{Observation}

\theoremstyle{definition}
\newtheorem{problem}[algorithm]{Problem}
\newtheorem{conj}[algorithm]{Conjecture}
\newtheorem{defin}[algorithm]{Definition}

\newtheorem{example}[algorithm]{Example}

\newcounter{stepenv}
\newenvironment{stepenv}[1][]{\refstepcounter{stepenv}}{}

\newcounter{step}[stepenv]

\newcounter{substep}[step]
\renewcommand{\thesubstep}{\thestep.\arabic{substep}}

\newcounter{claim}[stepenv]

\newcommand{\cF}{\mathcal{F}}

\newcommand{\cH}{\mathcal{H}}

\newcommand{\cO}{\mathcal{O}}

\newcommand{\cT}{\mathcal{T}}

\newcommand{\cW}{\mathcal{W}}

\newcommand{\bN}{\mathbb{N}}

\def\eps{{\varepsilon}}

\def\sm{\setminus}
\newcommand{\Set}[1]{\{#1\}}
\newcommand{\set}[2]{\{#1\,:\;#2\}}
\def\In{\subseteq}

\def\COMMENT#1{}
\def\TASK#1{}
\let\TASK=\footnote             

\begin{document}
	
	\title{Extremal aspects of graph and hypergraph decomposition problems}
	
	\date{\today}
	
	\author[S.~Glock]{Stefan Glock}
	\address[S.~Glock]{Institute for Theoretical Studies, ETH Z\"urich, CH}
	\email{dr.stefan.glock@gmail.com}
	
	\author[D.~K\"uhn]{Daniela K\"uhn}
	\address[D.~K\"uhn]{School of Mathematics, University of Birmingham, UK}
	\email{d.kuhn@bham.ac.uk}

	\author[D.~Osthus]{Deryk Osthus}
	\address[D.~Osthus]{School of Mathematics, University of Birmingham, UK}
	\email{d.osthus@bham.ac.uk}

	%

	\begin{abstract} 
		\noindent
		We survey recent advances in the theory of graph and hypergraph decompositions, with a focus on extremal results involving minimum degree conditions. We also collect a number of intriguing open problems, and formulate new ones.
	\end{abstract}
	
	\maketitle

\section{Introduction}\label{sec:intro}
The problem of decomposing large objects into (simple) smaller ones 
pervades many areas of Mathematics.
This is particularly true for Combinatorics.
Here, we focus on graphs and hypergraphs: given two graphs $F$ and $G$, an \emph{$F$-decomposition of $G$} is a collection of subgraphs of $G$, each isomorphic to $F$, such that every edge of $G$ is used exactly once.

Graph decomposition problems have a long history, dating back to the work of Euler on orthogonal Latin squares.
In 1847, Kirkman~\cite{kirkman:47} proved that the complete graph $K_n$ has a $K_3$-decomposition if and only if $n\equiv 1,3\mod{6}$. 
A related problem posed by Kirkman asks for a $K_3$-decomposition of $K_n$ such that the triangles can be organised into edge-disjoint $K_3$-factors, where a \emph{$K_3$-factor} is a set of vertex-disjoint triangles covering all vertices.
Also in the 19th century, Walecki proved the existence of decompositions of the complete graph $K_n$ into edge-disjoint Hamilton cycles (for odd $n$) and Hamilton paths (for even $n$). Note here that the assumptions on the parity of $n$ are necessary. Indeed, $K_n$ has $n(n-1)/2$ edges, and a Hamilton cycle has $n$ edges, so a decomposition into edge-disjoint Hamilton cycles has to consist of $(n-1)/2$ such cycles, implying that $n$ is odd. Similar `necessary divisibility conditions' can be observed for essentially every decomposition problem, and we will encounter many of these throughout this survey.

Classical graph decomposition results have mostly been obtained based on the symmetry
of the underlying structures, thus often involving algebraic techniques.
Recently, much progress has been made in the area of decompositions using probabilistic techniques.
This has gone hand in hand with the realisation that for many questions, it is not necessary that the underlying structure to be decomposed is `complete' or 
`highly symmetric'. 
This leads to the consideration of extremal aspects of such questions:
\emph{Are there natural (density) conditions which ensure (subject to the divisibility conditions) the existence of such decompositions?}

For Hamilton cycles and perfect matchings, the above question was resolved in~\cite{CBKLOT:16}, proving the so-called
Hamilton decomposition conjecture and the $1$-factorization conjecture for large~$n$:
the former states that for $d \ge \lfloor n/2 \rfloor$, every $d$-regular $n$-vertex graph $G$ has a decomposition into Hamilton cycles and at most one perfect matching, the latter states that the corresponding threshold for decompositions into perfect matchings
is $d \ge 2 \lceil n/4 \rceil -1$.
Here a decomposition into perfect matchings is often called a \emph{$1$-factorization}.
Similarly, Kelly's conjecture on Hamilton decompositions of regular tournaments not only turned out to be correct, 
but such a decomposition already exists in regular oriented graphs of 
minimum semi-degree at least $(3/8+o(1))n$~\cite{KO:13}.
For triangle decompositions (i.e.~Steiner triple systems), the corresponding question translates to a conjecture by 
Nash-Williams, which we introduce at the end of this section.

Minimum degree versions of decomposition problems have a natural application to `completion problems':
For instance, the $1$-factorization conjecture has an interpretation in terms of scheduling round-robin tournaments
(where $n$ players play all of each other in $n-1$ rounds): one can schedule the
first half of the rounds arbitrarily before one needs to plan the remainder of the tournament. 
More generally, there are applications to completions of partial designs, to hypergraph Euler tours and to the completion of Latin squares, to mention only a few.
We will discuss some of these applications in this survey.

Another feature of minimum degree conditions is that they provide a large 
class of graphs (or hypergraphs) whose membership is easy to verify 
algorithmically. In general, the question of whether a given \mbox{(hyper-)graph} $G$
has some decomposition into a given class of subgraphs is NP-complete~\cite{DT:97}.
Thus it is natural to seek simple sufficient conditions for the existence of such 
decompositions. Indeed, the minimum degree requirement yields attractive
results and conjectures in this respect. Another very fruitful notion is that of
quasirandomness -- we will briefly discuss this here too (mainly in the guise of
`typicality'), but will omit a detailed discussion.

We will focus mainly on decompositions into small subgraphs.
(We will also discuss related topics such as decompositions into $2$-factors and spanning trees. For Hamilton decompositions, we refer e.g.~to~\cite{KO:14c}).
The main conjecture in this area is due to Nash-Williams, which we introduce now: 
Observe that if a graph $G$ admits a $K_3$-decomposition, then the number of edges of $G$ must be divisible by~$3$, and all the vertex degrees of $G$ must be even. We say that $G$ is \emph{$K_3$-divisible} if it has these properties. Clearly, not every $K_3$-divisible graph has a $K_3$-decomposition (e.g.~$C_6$). In fact, to decide whether a given graph has a $K_3$-decomposition is NP-hard~\cite{DT:97}.
However, the following beautiful conjecture of Nash-Williams suggests that if the minimum degree of $G$ is sufficiently large, then the existence of a $K_3$-decomposition hinges only on the necessary divisibility condition.

\begin{conj}[\cite{nash-williams:70}]  \label{conj:NW}
	For sufficiently large $n$, every $K_3$-divisible graph $G$ on $n$ vertices with $\delta(G)\ge 3n/4$ has a $K_3$-decomposition.
\end{conj}

The following class of extremal examples shows that the bound on the minimum degree would be best possible. 

\begin{figure}
	\begin{center}
	\begin{tikzpicture}[scale=0.8, line width=1pt]
	
	\begin{scope}[shift={(0,0)}]
	\filldraw[black!30] (-2,-0.6)--(-2,0.6)--(2,0.6)--(2,-0.6)--(-2,-0.6); 
	\end{scope}
	\begin{scope}[shift={(0,-4)}]
	\filldraw[black!30] (-2,-0.6)--(-2,0.6)--(2,0.6)--(2,-0.6)--(-2,-0.6); 
	\end{scope}
	\begin{scope}[shift={(-2,-2)},rotate=90]
	\filldraw[black!30] (-2,-0.6)--(-2,0.6)--(2,0.6)--(2,-0.6)--(-2,-0.6); 
	\end{scope}
	\begin{scope}[shift={(2,-2)},rotate=90]
	\filldraw[black!30] (-2,-0.6)--(-2,0.6)--(2,0.6)--(2,-0.6)--(-2,-0.6); 
	\end{scope}
	
	\filldraw[fill = white] (-2,0) circle  (0.8);
	\filldraw[fill = white] (2,0) circle  (0.8);
	\filldraw[fill = white] (-2,-4) circle  (0.8);
	\filldraw[fill = white] (2,-4) circle  (0.8);

	\filldraw[fill = black!40,pattern = crosshatch dots] (-2,0) circle  (0.8);
	\filldraw[fill = black!40,pattern= north west lines] (2,0) circle  (0.8);
	\filldraw[fill = black!40,pattern= north west lines] (-2,-4) circle  (0.8);
	\filldraw[fill = black!40,pattern= crosshatch dots] (2,-4) circle  (0.8);
	
	\end{tikzpicture}
\end{center}
	\caption{The standard example is a blown-up $C_4$, where $G_1$ and $G_2$ each consist of two disjoint cliques.}
	\label{fig:NW}
\end{figure}
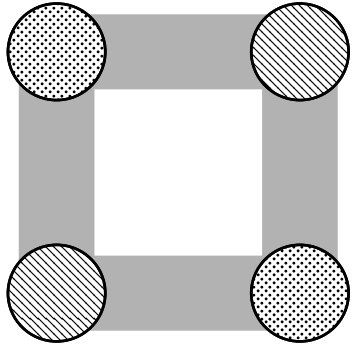

\begin{example}\label{ex:NW}
	Given any $k\in \bN$, let $G_1$ and $G_2$ be vertex-disjoint $(6k+2)$-regular graphs with $|G_1|=|G_2|=12k+6$ and let $G_3$ be the complete bipartite graph between $V(G_1)$ and $V(G_2)$. Let $G:=G_1\cup G_2\cup G_3$. (In the standard construction, each of $G_1$ and $G_2$ is a union of two disjoint cliques of size~$6k+3$, see Figure~\ref{fig:NW}.)
	It is straightforward to check that $\delta(G)=3|G|/4-1$\COMMENT{$|G|=24k+12$. $\delta(G)=6k+2+12k+6=18k+8$.} and that $G$ is $K_3$-divisible.\COMMENT{all degrees even, $e(G_3)=(12k+6)^2$, $e(G_1\cup G_2)=(12k+6)(6k+2)$} However, every triangle in $G$ contains at least one edge from $G_1\cup G_2$. Since $2e(G_1\cup G_2)<e(G_3)$, $G$ cannot have a $K_3$-decomposition.\footnote{This type of extremal example is called a `space barrier'. There are simply not enough edges in $G_1\cup G_2$ for a triangle decomposition. In other constructions, the obstacle might be a `divisibility barrier'.}
\end{example}

Conjecture~\ref{conj:NW} is still open. Nevertheless, it has already inspired very fruitful research, and will serve us as a thread running through this survey.

\subsection{Organisation of this survey}

In Section~\ref{sec:notation}, we collect some basic notation.
In Section~\ref{sec:fractional}, we consider approximate versions of Conjecture~\ref{conj:NW}  and its generalisations, and in Section~\ref{sec:thresholds} we discuss how such approximate results can be turned into exact ones.
In Section~\ref{sec:Fdecs}, we introduce the decomposition problem for hypergraphs. 
Subsequently, in Section~\ref{sec:euler} we present an application of hypergraph decompositions to Euler tours.
Section~\ref{sec:OW} is devoted to the Oberwolfach problem.
Finally, in Section~\ref{sec:related}, we close by briefly mentioning some further open problems.

\subsection{Notation}\label{sec:notation}
\COMMENT{Always use $k$ for the uniformity if possible, and $r$ for clique size. Use $d$ for regularity of host graphs.}

Let us briefly agree on some general notation.
A \emph{hypergraph} $G$ consists of a set of vertices $V(G)$ and a set of edges $E(G)$, where each edge is a subset of the vertex set. The hypergraph is called \emph{$k$-uniform}, or just a \emph{$k$-graph}, if every edge has size~$k$. Hence, a $2$-graph is simply a graph.
We let $|G|$ denote the number of vertices of $G$ and $e(G)$ the number of edges.

Let $G$ be a $k$-graph.
For a set $S\In V(G)$ with $|S| = k-1$, we let $N_G(S)$ denote the \emph{neighbourhood of $S$ in $G$}, that is, the set of vertices which form an edge
in $G$ together with~$S$. For a set $S\In V(G)$ with $0 \le |S| \le k-1$, 
the \emph{degree $d_G(S)$ of $S$} is the number of edges in $G$
containing~$S$. 
We let $\delta(G)$ and $\Delta(G)$ denote the minimum and maximum $(k-1)$-degree of $G$, respectively, that is, the minimum/maximum value of $d_G(S)$ over all $S\In V(G)$ of size~$k-1$.

Some of the results we state apply for `quasirandom' (hyper-)graphs, by which we mean the following notion.

\begin{defin}[Typicality]\label{def:typical}
	A $k$-graph $G$ on $n$ vertices is called \emph{$(c,h,p)$-typical} if for any set $A$ of $(k-1)$-subsets of $V(G)$ with $|A|\le h$, we have $|\bigcap_{S\in A}N_G(S)| = (1\pm c)p^{|A|}n$. 
\end{defin}

We let $K^k_n$ denote the complete $k$-graph on $n$ vertices.
Finally, $[n]$ denotes the set $\Set{1,\dots,n}$.

\section{Approximate and fractional decompositions} \label{sec:fractional}

As Conjecture~\ref{conj:NW} turned out to be hard, one might try to prove at least an approximate version.
For instance, one could attempt to obtain a collection of edge-disjoint triangles which do not cover all edges of $G$, but may leave $o(|G|^2)$ edges uncovered. We refer to such types of `almost' decompositions as \emph{approximate decompositions}.
Another relaxation is to find only a `fractional' decomposition.
As we shall see, these two concepts are closely related, and played a pivotal role in recent progress towards Conjecture~\ref{conj:NW} and many related problems.

\subsection{From fractional to approximate decompositions}

Fractional relaxations of graph parameters have been studied extensively in recent decades. Often, it turns out that these fractional relaxations give good approximations for the original parameter. This can also have important algorithmic ramifications. 

A \emph{fractional $F$-decomposition} of a graph $G$ is a function $\omega$ that assigns to each copy of $F$ in $G$ a value in $[0,1]$ such that for all $e \in E(G)$,
\begin{align*} 
\sum_{F'\colon e\in E(F')}  \omega(F') = 1,
\end{align*}
where the sum is over all copies $F'$ of $F$ which contain~$e$.

Thus, an $F$-decomposition is a fractional $F$-decomposition taking values in~$\Set{0,1}$. Note that it can be much easier to obtain a fractional $F$-decomposition rather than an $F$-decomposition. For example, $K_n$ always has a fractional $F$-decomposition (assuming $n\ge |F|$ of course), by giving every copy of $F$ the same weight. 

The following theorem of Haxell and R\"odl~\cite{HR:01} allows us to turn a fractional decomposition into an approximate decomposition.

\begin{theorem}[\cite{HR:01}]\label{thm:HR}
	If an $n$-vertex graph $G$ has a fractional $F$-decomposition, then all but $o(n^2)$ edges of $G$ can be covered with edge-disjoint copies of~$F$.
\end{theorem}

The proof is based on Szemer\'edi's regularity lemma and a theorem of Frankl and R\"odl~\cite{FR:85}. The latter result allows one to obtain an approximate $F$-decomposition of a dense graph $G$ whenever every edge of $G$ is contained in roughly the same number of copies of~$F$. This is a special case of a more general theorem on almost perfect matchings in (almost) regular hypergraphs with small maximum codegree, proved using the celebrated `R\"odl nibble'.
Roughly speaking, the main idea in the proof of Theorem~\ref{thm:HR} is to partition $G$ into edge-disjoint subgraphs, each of which has the above regularity condition. Such a partition is obtained with the help of Szemer\'edi's regularity lemma, using the assumption that $G$ has a fractional $F$-decomposition.
One can then apply the Frankl--R\"odl result to all these subgraphs individually to obtain the desired approximate $F$-decomposition.

Given that there are numerous algorithmic applications,
it would be very interesting to obtain a proof of Theorem~\ref{thm:HR} which 
does not rely on Szemer\'edi's regularity lemma.

Theorem~\ref{thm:HR} has been generalised to hypergraphs by R\"{o}dl, Schacht, Siggers and Tokushige~\cite{RSST:07}.

\subsection{Fractional decomposition thresholds}\label{subsec:fracs}

Motivated by this, we could aim to prove the fractional version of Conjecture~\ref{conj:NW}. 
Observe that Example~\ref{ex:NW} also shows that the same minimum degree condition would be optimal for this relaxed version.

We define the \emph{fractional $F$-decomposition threshold $\delta_F^\ast$} to be the infimum of all $\delta\in[0,1]$ with the following property: there exists $n_0\in \bN$ such that any $F$-divisible graph $G$ on $n\ge n_0$ vertices with $\delta(G)\ge \delta n$ has a fractional $F$-decomposition.
We remark that it is somehow unnatural to assume divisibility here since a graph can have a fractional $F$-decomposition without being $F$-divisible. The sole purpose of this assumption is to have the inequality $\delta_F^*\le \delta_F$ hold trivially, where $\delta_F$ is the threshold for $F$-decompositions which we will define in Section~\ref{sec:thresholds}. This underpins our motivation that fractional decompositions are a relaxation of the original decomposition problem.\footnote{Having said this, we point out that all the results in this section on $\delta_F^*$ also hold without the divisibility assumption, that is, none of the proofs showing existence of fractional $F$-decompositions make any use of the divisibility assumption, and we suspect that the value of $\delta_F^*$ is the same for both variants.}

Until recently, the best bound on $\delta_{K_3}^\ast$ was obtained by Dross~\cite{dross:16},
who showed that $\delta_{K_3}^\ast \le 0.9$, using an elegant approach based on the
max-flow-min-cut theorem. This improved earlier work of Yuster~\cite{yuster:05}, as well as
Dukes~\cite{dukes:12,dukes:15} and Garaschuk~\cite{garaschuk:14}.
Delcourt and Postle~\cite{DP:19} hold the current record (a slightly weaker bound was obtained simultaneously by Dukes and Horsley~\cite{DH:20}).
An obvious open problem is to obtain further significant improvements.

\begin{theorem}[\cite{DP:19}]
	$\delta_{K_3}^\ast \le (7+\sqrt{21})/14 \approx 0.82733$.
\end{theorem}

The best current bound on the fractional decomposition threshold of larger cliques was proved by Montgomery~\cite{montgomery:19}.

\begin{theorem}[\cite{montgomery:19}]
	For each $r\ge 4$, we have $\delta_{K_r}^\ast \le 1-1/(100r)$.
\end{theorem}

This result is best possible up to the constant~$100$. Indeed, a natural generalisation of Example~\ref{ex:NW} (see~\cite{yuster:05}) shows that the following conjecture (which is implicit in~\cite{gustavsson:91}) would be optimal.

\begin{conj}\label{conj:frac clique}
	For each $r\ge 3$, we have $\delta_{K_r}^\ast \le 1-1/(r+1)$.
\end{conj}

So far, we have only stated bounds for the fractional decomposition threshold of cliques. Yuster~\cite{yuster:12} showed that these bounds can be used to obtain approximate $F$-decompositions for general $r$-chromatic graphs~$F$.

\begin{theorem}[\cite{yuster:12}]\label{thm:clique2F}
	For any fixed graph $F$ with $r=\chi(F)$, any sufficiently large $n$-vertex graph $G$ with $\delta(G)\ge ( \delta_{K_r}^*+o(1))n$ contains edge-disjoint copies of~$F$ covering all but $o(n^2)$ edges of $G$.
\end{theorem}

\subsection{Bandwidth theorem for approximate decompositions}

Not only does $\delta_{K_r}^*$ yield minimum degree bounds for approximate decompositions into a fixed $r$-chromatic graph $F$, but even for spanning $r$-chromatic graphs~$F$, provided that $F$ has bounded degree, the chromatic number of $F$ is `essentially' equal to $r$ and $F$ is `path-like'. More precisely,
a graph $F$ is \emph{$(r,\eta)$-chromatic} if the graph $F'$, obtained from $F$ by deleting isolated vertices, can be properly coloured with $r+1$ colours such that one colour class has size at most $\eta |F'|$.
For instance, the cycle $C_\ell$ is $(2,1/\ell)$-chromatic.
Moreover, $F$ is \emph{$\eta$-separable} if there exists a set $S$ of at most $\eta |F|$ vertices such that each component of $F\sm S$ has size at most $\eta |F|$. For bounded degree graphs, being separable is equivalent to having sublinear bandwidth.
Examples of separable graphs include cycles, powers of cycles, trees and planar graphs.

The following result by
Condon, Kim, K\"uhn and Osthus~\cite{CKKO:19} provides a degree condition which ensures that a regular graph $G$ has an approximate
decomposition into $\cF$ for any collection $\cF$ of $(r,\eta)$-chromatic $\eta$-separable graphs of bounded degree.
The degree condition is best possible in general (unless one has additional information about
the graphs in $\cF$). The `original' bandwidth theorem of B\"ottcher, Schacht and Taraz~\cite{BST:09} involves a similar degree condition for embedding a single
spanning subgraph.

Given a collection $\cF$ of graphs, we let $e(\cF)=\sum_{F\in \cF}e(F)$. Moreover, we say that $\cF$ \emph{packs} into a graph $G$ if the graphs from $\cF$ can be embedded edge-disjointly into~$G$.

\begin{theorem}[\cite{CKKO:19}]\label{thm:bandwidth}
	For all $\Delta,r\ge 2$ and $\eps>0$, there exist $\eta>0$ and $n_0\in \bN$ such that the following holds for all $n\ge n_0$. Assume that $\cF$ is a collection of $(r,\eta)$-chromatic, $\eta$-separable $n$-vertex graphs with maximum degree at most~$\Delta$.
	Assume that $G$ is a $d$-regular $n$-vertex graph with $d\ge (\max\Set{\delta_{K_r}^*,1/2} + \eps) n$.
	If $e(\cF)\le (1-\eps)e(G)$, then $\cF$ packs into~$G$.
\end{theorem}

The proof of Theorem~\ref{thm:bandwidth} is based on the blow-up lemma for approximate decompositions developed in~\cite{KKOT:19}, a subsequent shorter proof of which can be found in~\cite{EJ:20}.

Since $\delta_{K_2}^*=0$ and trees are separable, Theorem~\ref{thm:bandwidth} has the following immediate consequence for approximate decompositions into bounded degree trees.

\begin{cor}[\cite{CKKO:19}] \label{cor:treebandwidth}
	For all $\eps,\Delta>0$, the following holds for sufficiently large~$n$. Let $\cT$ be a collection of trees on at most $n$ vertices with maximum degree at most~$\Delta$.
	Assume that $G$ is a $d$-regular $n$-vertex graph with $d\ge (1/2+\eps) n$. If $e(\cT)\le (1-\eps)e(G)$, then $\cT$ packs into~$G$.
\end{cor}

Similarly, $2$-regular graphs are separable and have chromatic number at most~$3$. Recall that $\delta_{K_3}^*$ is conjectured to be $3/4$, though the currently best known upper bound is roughly~$0.83$.

\begin{cor}[\cite{CKKO:19}] \label{cor:OW mindeg approx}
	For all $\eps>0$, the following holds for sufficiently large~$n$. Assume that $\cF$ is a collection of $2$-regular $n$-vertex graphs.
	Assume that $G$ is a $d$-regular $n$-vertex graph with $d\ge (\delta_{K_3}^* + \eps) n$. If $e(\cF)\le (1-\eps)e(G)$, then $\cF$ packs into~$G$.
\end{cor}
This corollary was instrumental in the resolution of the Oberwolfach problem
(see Section~\ref{sec:OW}).

Moreover, if a $2$-regular graph has only long cycles, say of length at least $\ell$, then it is $(2,1/\ell)$-chromatic.

\begin{cor}[\cite{CKKO:19}]\label{cor:OW mindeg approx girth}
	For all $\eps>0$, there exist $\ell,n_0\in \bN$ such that the following holds for all $n\ge n_0$. Assume that $\cF$ is a collection of $2$-regular $n$-vertex graphs with girth at least~$\ell$.
	Assume that $G$ is a $d$-regular $n$-vertex graph with $d\ge (1/2+\eps) n$. If $e(\cF)\le (1-\eps)e(G)$, then $\cF$ packs into~$G$.
\end{cor}

Later in Section~\ref{sec:OW}, we will discuss possible `exact' versions (i.e.~`full' decompositions)
of Corollaries~\ref{cor:OW mindeg approx} and~\ref{cor:OW mindeg approx girth} in the context of the Oberwolfach problem.

\section{Decomposition thresholds for fixed graphs \texorpdfstring{$F$}{F}}\label{sec:thresholds}

Recall that a graph $G$ is $K_3$-divisible if $3\mid e(G)$ and $2\mid d_G(v)$ for all $v\in V(G)$, and that any graph with a $K_3$-decomposition must be $K_3$-divisible.
Now, we generalise this to arbitrary graphs.
Let $F$ be a fixed graph. Define $\gcd(F)$ as the greatest common divisor of all vertex degrees of~$F$.
We say that a graph $G$ is \emph{$F$-divisible} if $e(F)\mid e(G)$ and $\gcd(F)\mid \gcd(G)$. It is easy to see that a necessary condition for $G$ to have an $F$-decomposition is to be $F$-divisible.

In 1976, Wilson~\cite{wilson:75,wilson:76} proved the following fundamental result.
\begin{theorem}[\cite{wilson:75,wilson:76}]\label{thm:wilson}
	Given any graph~$F$, for sufficiently large~$n$, the complete graph $K_n$ has an $F$-decomposition whenever it is $F$-divisible.
\end{theorem}
Gustavsson~\cite{gustavsson:91} showed that the same is true with $K_n$ replaced by an `almost-complete' graph~$G$, where $\delta(G)\ge (1-\eps)|G|$ for some tiny $\eps>0$.\footnote{This proof has not been without criticism.}

We define the \emph{$F$-decomposition threshold $\delta_F$} to be the infimum of all $\delta\in[0,1]$ with the following property: there exists $n_0\in \bN$ such that any $F$-divisible graph $G$ on $n\ge n_0$ vertices with $\delta(G)\ge \delta n$ has an $F$-decomposition.

Thus, Conjecture~\ref{conj:NW} would imply that $\delta_{K_3}\le 3/4$ (which would be tight by Example~\ref{ex:NW}).
Of course, it is interesting to investigate the decomposition threshold of arbitrary graphs, not just triangles.

\begin{problem}
	Determine the $F$-decomposition threshold for every graph~$F$.
\end{problem}

The remainder of this section is devoted to this problem. 
In~\cite{BKLO:16}, a general method was developed that turns an approximate decomposition into an exact decomposition. 
This approach was refined in \cite{GKLMO:19} to yield the following result, which gives a general upper bound on the decomposition threshold. 

\begin{theorem}[\cite{GKLMO:19}]\label{thm:dec main upper}
	For any graph $F$, we have $$\delta_F\le \max\Set{\delta_F^*,1-1/(\chi(F)+1)}.$$
\end{theorem}

Theorem~\ref{thm:dec main upper} improves a bound of $\delta_F\le \max\Set{\delta_F^*,1-1/3r}$ proved in \cite{BKLO:16} for $r$-regular graphs~$F$. Also, the cases where $F=K_3$ or $C_4$ were already proved in~\cite{BKLO:16}.

Since $\delta_{K_r}^* \ge 1-1/(r+1)$ (see Section~\ref{subsec:fracs}),
Theorem~\ref{thm:dec main upper} implies that the decomposition threshold for cliques equals its fractional relaxation.

\begin{cor}[\cite{GKLMO:19}] \label{cor:equal thresholds}
	For all $r\ge 3$, $\delta_{K_r}=\delta_{K_r}^\ast$.
\end{cor}

An affirmative answer to Conjecture~\ref{conj:frac clique} would, via (a variant of) Theorem~\ref{thm:dec main upper} and Theorem~\ref{thm:clique2F}, also imply the following general upper bound on the $F$-decomposition threshold.

\begin{conj}[\cite{GKLMO:19}]
	For every graph $F$, we have $\delta_F\le 1-1/(\chi(F)+1)$.
\end{conj}

This would be tight in general. However, a natural question is whether the bound can be improved for certain graphs~$F$. This will be the topic of Sections~\ref{sec:bip dec} and~\ref{sec:discrete}.
Before, we discuss very roughly the idea behind the proof of Theorem~\ref{thm:dec main upper}.

\subsection{Turning approximate decompositions into exact ones}

Let $G$ be a large $F$-divisible graph. We always assume the minimum degree to be at least $(\delta_F^* +o(1))|G|$. By definition of $\delta_F^*$, this gives us a fractional $F$-decomposition of $G$ for free, and Theorem~\ref{thm:HR} then provides us with an approximate $F$-decomposition. The challenge is to deal with the leftover of such an approximate decomposition. For this, the following concept is crucial.
\begin{defin}\label{def:absorber}
	An \emph{$F$-absorber for a graph $L$} is a graph $A$ which is edge-disjoint from $L$ such that both $A$ and $A\cup L$ have an $F$-decomposition.
\end{defin}
(Here, one can think of both $A$ and $L$ being subgraphs of a given host graph, but this host graph plays no role in the definition.)
This motivates the following naive strategy: For any `possible leftover' $L$, find an $F$-absorber $A_L$ in $G$ such that all these absorbers are edge-disjoint. Let $A$ be the union of these `exclusive' absorbers.
Then obtain an approximate decomposition of $G-E(A)$. Now, we are guaranteed that the leftover will be one of the graphs $L$ for which we found an absorber~$A_L$. We can use the fact that $A_L\cup L$ has an $F$-decomposition by Definition~\ref{def:absorber}. Similarly, all other absorbers also have an $F$-decomposition themselves, which gives in total an $F$-decomposition of~$G$.

The astute reader will already have noticed (at least) one major problem with this approach. The number of possible leftovers $L$ is gigantic. Although we know that $e(L)=o(|G|^2)$, there are still exponentially many possibilities, and so there is no hope to find edge-disjoint(!) absorbers for each of them. 

However, by using an `iterative absorption' process, one can successively obtain approximate decompositions of smaller and smaller pieces of $G$, combined with some partial absorption steps, until finally the number of possible leftovers is so small that the above naive strategy works as the last step. We do not go into more detail for this.
Ultimately, the remaining problem one has to solve is the following: Given $G$ and a small subgraph $L\In G$ (we can even assume that $|L|$ is bounded), find an $F$-absorber for $L$ in~$G$. This can be achieved if $\delta(G)\ge (1-1/(\chi(F)+1) +o(1))|G|$, which together with our initial assumption (needed for approximate decompositions) yields the bound in Theorem~\ref{thm:dec main upper}.

The expository article \cite{BGKLMO:20} contains a short proof, with all details, of Theorem~\ref{thm:dec main upper} in the case when $F$ is a triangle.

\subsection{Bipartite graphs}\label{sec:bip dec}

One big advantage when dealing with bipartite graphs is that we do not have to worry about approximate decompositions. Indeed, if $F$ is bipartite, then its Tur\'an density vanishes, which means that we can iteratively pull out copies of $F$ from a graph $G$ until at most $o(|G|^2)$ edges remain.

In~\cite{GKLMO:19}, the decomposition threshold $\delta_F$ was determined for every bipartite graph~$F$. Before stating the general result, let us mention preceding work.
Yuster~\cite{yuster:98} showed that $\delta_F=1/2$ if $F$ is a tree, and later generalised his own result by drawing the same conclusion if $F$ is connected and contains a vertex of degree one~\cite{yuster:02}.
Bryant and Cavenagh~\cite{BC:15} proved that $\delta_{C_4}\le 31/32$.
Barber, K\"uhn, Lo and Osthus~\cite{BKLO:16} showed that $\delta_{K_{r,r}} \le 1-1/(r+1)$ and $\delta_{C_\ell}=1/2$ for all even $\ell\ge 6$. Note that their first result implies $\delta_{C_4} \le 2/3$. The following extremal example due to Taylor (see~\cite{BKLO:16}) shows that this is tight, giving $C_4$ a special role among the even cycles.

\begin{example}\label{ex:C4dec}
	Consider the following graph $G$: Let $A,B,C$ be sets of size roughly $n/3$ where $G[A],G[C]$ are complete, $B$ is independent and $G[A,B]$ and $G[B,C]$ are complete bipartite. Clearly, $\delta(G)\approx 2n/3$. It is easy to see that any $C_4$ in $G$ contains an even number of edges from~$A$. Thus, if $e_G(A)$ is odd, then $G$ cannot have a $C_4$-decomposition. By slightly altering the sizes of $A,B,C$, it is not difficult to ensure this while also ensuring that $G$ is $C_4$-divisible. 
\end{example}

From Theorem~\ref{thm:dec main upper}, we have that $\delta_F\le 2/3$ for all bipartite graphs~$F$.\footnote{Note that Theorem~\ref{thm:dec main upper} is formulated in terms of the fractional decomposition threshold. However, for bipartite graphs, the detour via fractional decompositions is not necessary as explained at the beginning of this section.} We have seen that this is tight for $C_4$, so a natural question is whether we can show this for other graphs too.
Example~\ref{ex:C4dec} motivates the following definition. A set $X\In V(F)$ is called \emph{$C_4$-supporting in $F$} if there exist distinct $a,b\in X$ and $c,d\in V(F)\sm X$ such that $ac,bd,cd\in E(F)$. Observe that if $F$ is a copy of $C_4$ in $G$ as in Example~\ref{ex:C4dec}, then $V(F)\cap A$ is not $C_4$-supporting in~$F$.
We define
\begin{align*}
\tau(F)&:=\gcd\set{e(F[X])}{X\In V(F)\text{ is not }C_4\text{-supporting in }F}.
\end{align*}
For instance, $\tau(C_4)=2$. Whenever $\tau(F)>1$, we can adapt Example~\ref{ex:C4dec} as follows: ensure that $e_G(A)$ is not divisible by $\tau(F)$, then $G$ cannot have an $F$-decomposition by the same logic. At the same time, ensure that $G$ is $F$-divisible. (This needs some additional adjustment, but is a mere technicality, so we omit the details.)
We note that $\tau(F)\mid \gcd(F)$ (consider all sets $X$ consisting of a vertex and its neighbourhood).

If $\tau(F)=1$, which happens for instance if there exists an edge in $F$ that is not contained in any~$C_4$, or if $\gcd(F)=1$, then the construction fails.
In fact, it turns out that in this case, the special structure of $F$ can be exploited to design absorbers which show that $\delta_F\le 1/2$. This is tight if $F$ is connected and $e(F)\ge 2$. 
In order to deal with disconnected graphs, we further define 
\begin{align*}
\tilde{\tau}(F)&:=\gcd\set{e(C)}{C\text{ is a component of }F}.
\end{align*}
If $\tilde{\tau}(F)>1$, we can use as an extremal example the disjoint union of two cliques on roughly $n/2$ vertices (modulo some simple modifications to ensure divisibility). Moreover, if $\tilde{\tau}(F)=1$, by adding one `cross-edge' between the cliques, it is clear that there cannot be an $F$-decomposition if every edge of $F$ is contained in a cycle.

Altogether, we have the following complete picture.

\begin{theorem}[\cite{GKLMO:19}] \label{thm:bipartite char}
	Let $F$ be a bipartite graph. Then
	$$\delta_F = \begin{cases}
	2/3 &\mbox{if } \tau(F)>1; \\
	0 &\mbox{if } \tilde{\tau}(F)=1 \mbox{ and } F \mbox{ has a bridge}; \\
	1/2 &\mbox{otherwise. }
	\end{cases}$$
\end{theorem}

In particular, this implies that if $F$ is a bipartite connected graph with at least two edges, then $\delta_F=1/2$ if $\gcd(F)=1$ or $F$ has an edge that is not contained in any~$C_4$.
Moreover, since $\tau(K_{s,t})=\gcd(s,t)$, we have $\delta_{K_{s,t}}=1/2$ if $s$ and $t$ are coprime and $\delta_{K_{s,t}}=2/3$ otherwise.

Given that Theorem~\ref{thm:bipartite char} asymptotically determines the decomposition
threshold for bipartite graphs, the next question is of course whether one can determine
the exact threshold. This was achieved for even cycles (with the exception of 
$6$-cycles) by Taylor~\cite{taylor:19} and for trees 
(for infinitely many $n$) by Yuster~\cite{yuster:00b}.
It would be particularly interesting to solve the cases when $F$ is complete bipartite
and when $F$ is a $6$-cycle.

\subsection{A discretisation result}\label{sec:discrete}

Recall that $\delta_F\le \max\Set{\delta_F^*,1-1/(\chi(F)+1)}$ holds, and probably $\delta_F\le 1-1/(\chi(F)+1)$ for every graph~$F$. In the previous subsection, we saw that if $\chi:=\chi(F)=2$, then $\delta_F\in \Set{1-1/(\chi-1),1-1/\chi,1-1/(\chi+1)}$. The question suggests itself whether a similar phenomenon occurs in general. (Note that we have $\delta_F\ge 1-1/(\chi-1)$ as shown by the Tur\'an graph.) 
Progress towards this has been made in~\cite{GKLMO:19} for graphs of chromatic number at least~$5$.

\begin{theorem}[\cite{GKLMO:19}]\label{thm:dec discrete}	
	Let $F$ be a graph with $\chi:=\chi(F)\ge 5$. Then we have $\delta_F\in\Set{\delta_F^\ast,1-1/\chi,1-1/(\chi+1)}$.
\end{theorem}

The proof of this result has the nice feature that the assumption $\delta_F < 1-1/(\chi+1)$ is used indirectly to obtain better absorbers. 
Very roughly speaking, the idea is as follows:
Suppose that $\delta_F < 1-1/(\chi+1)$. By the definition of $\delta_F$, this means that a graph which is very close to a large balanced complete $(\chi+1)$-partite  graph has an $F$-decomposition if it is $F$-divisible. One can use such graphs as absorbers, or more precisely, as building blocks of absorbers. Although they might be very large in size, the fact that they are $(\chi+1)$-partite allows us to find them in a given graph $G$ already if $\delta(G) \ge (1-1/\chi +o(1))|G|$, which we can use (modulo other parts of the argument) to show that $\delta_F\le 1-1/\chi$. Similarly, assuming $\delta_F < 1-1/\chi$ allows one to find absorbers if $\delta(G) \ge (1-1/(\chi-1) +o(1))|G|$. 
The complete proof of Theorem~\ref{thm:dec discrete}  involves a number of reductions for which the assumption $\chi\ge 5$ is needed.

It would be interesting to investigate the cases $\chi(F)\in \Set{3,4}$ further. Perhaps a good starting point are odd cycles.

\begin{problem}
	Determine $\delta_{C_{2\ell+1}}$. 
\end{problem} 

It is shown in~\cite[Proposition~12.1]{BKLO:16}  that $\delta_{C_{2\ell+1}}\ge \frac{1}{2} + \frac{1}{4\ell}$ (and the same construction gives the same bound for $\delta^*_{C_{2\ell+1}}$).
In the first version of this survey, we posed the problem of at least showing that $\delta_{C_{2\ell+1}} \to 1/2$ as $\ell \to \infty$. By the results in~\cite{BKLO:16}, it suffices to prove this for the fractional threshold. This was confirmed very recently by Joos and K\"uhn~\cite{JK:21}. Moreover, they obtain an analogous result for fractional decompositions of hypergraphs.

\subsection{Decompositions of partite graphs and Latin squares}

A \emph{Latin square} of order $n$ is an $n \times n$ array of cells, each containing a symbol from $[n]$, where no symbol appears twice in any row or column.
A Latin square  corresponds to a $K_3$-decomposition of the complete tripartite graph $K_{n,n,n}$ with vertex classes consisting of the rows, columns and symbols.

More generally, two Latin squares $R$ (red) and $B$ (blue) drawn in the same $n \times n$ array of cells are \emph{orthogonal} if no two cells contain the same combination of red symbol and blue symbol.
As above, it is easy to see that a pair of orthogonal Latin squares corresponds to a $K_4$-decomposition of $K_{n,n,n,n}$. 
Even more generally, there is a bijection between sequences of $r-2$ \emph{mutually orthogonal} 
Latin squares (where every pair from the sequence are orthogonal) and $K_r$-decompositions of balanced complete $r$-partite graphs.

A partial Latin square is defined similarly as a Latin square, except that cells are allowed to be empty.
Daykin and H\"aggkvist~\cite{DH:84} made the following conjecture:
\begin{conj}[\cite{DH:84}]\label{conj:LS complete}
	Given a partial Latin square $L$ of order $n$ where each row, column and symbol is used at most $n/4$ times, it is possible to complete $L$ into a Latin square of order~$n$. 
\end{conj}
This can be seen as a natural analogue of the conjecture
of Nash-Williams (Conjecture~\ref{conj:NW}).
Indeed, it turns out that they are closely related:
Barber, K\"uhn, Lo, Osthus and Taylor~\cite{BKLOT:17} proved an analogue of Corollary~\ref{cor:equal thresholds} for $r$-partite graphs.
Complemented by results of Bowditch and Dukes~\cite{BD:19} on fractional $K_3$-decompositions of
balanced tripartite graphs and results of Montgomery~\cite{montgomery:17} for larger cliques, this implies that a (sequence of mutually orthogonal) partial Latin squares can be completed provided that no row, column or coloured symbol has already been used too often (see~\cite{BKLOT:17} for a precise statement).

In particular, the case $r=3$ shows that Conjecture~\ref{conj:LS complete} holds with $1/25-o(1)$ instead of~$1/4$. 
\begin{theorem}[\cite{BKLOT:17,BD:19}]\label{thm:mols}
	Given a partial Latin square $L$ of order $n$ where each row, column and symbol is used at most $(\frac{1}{25}-o(1)) n$ times, it is possible to complete $L$ into a Latin square of order~$n$. 
\end{theorem}

It would be very interesting to improve these bounds
(recall that the obvious `bottleneck' consists of finding improved bounds for the 
fractional decomposition problem).

\section{\texorpdfstring{$F$}{F}-decompositions of hypergraphs}\label{sec:Fdecs}

So far, we have only considered decomposition problems for graphs, but of course the same type of questions can also be asked for hypergraphs. In fact, most definitions and questions translate almost verbatim to hypergraphs. However, the known results are much more limited.

Let $G$ and $F$ be $k$-graphs. An \emph{$F$-decomposition of $G$} is a collection of copies of $F$ in $G$ such that every edge of $G$ is contained in exactly one of these copies.
As for graphs, the existence of an $F$-decomposition necessitates certain divisibility conditions. For instance, we surely need $e(F)\mid e(G)$.
More generally, define $$d_F(i):=\gcd\set{d_F(S)}{S\In V(F),|S|=i}$$ for all $0\le i\le k-1$.\COMMENT{As long as $F$ is not edgeless, this is well-defined.} Note that $d_F(0)=e(F)$. Moreover, $d_F(1)=\gcd(F)$ for $k=2$.
Now, $G$ is called \emph{$F$-divisible} if $d_F(i)\mid d_G(i)$ for all $0\le i\le k-1$.
It is easy to see that $G$ must be $F$-divisible in order to admit an $F$-decomposition. 

In 1853, Steiner asked for which $n>r>k$ there exists a $K^{k}_r$-decomposition of~$K^k_n$. Clearly, the above divisibility conditions need to be satisfied. The folklore `Existence conjecture' postulated that these conditions are also sufficient, at least when $r$ and $k$ are fixed and $n$ is sufficiently large.
This was proved in a recent breakthrough by Keevash~\cite{keevash:14}.

\begin{theorem}[\cite{keevash:14}]
	For fixed $r,k$ and sufficiently large $n$, $K^k_n$ has a $K^{k}_r$-decomposition if it is $K^{k}_r$-divisible.
\end{theorem}

In fact, Keevash proved a more general result that holds also if $K^k_n$ is replaced by a typical $k$-graph~$G$.

An obvious question (e.g.~asked by Keevash) is whether $K^{k}_r$ can be replaced by any $k$-graph~$F$. 
For instance, Hanani~\cite{hanani:79} settled the problem if $F$ is an octahedron (viewed as a $3$-uniform hypergraph).

Recently, the authors together with Lo solved the general problem~\cite{GKLO:ta}, thus extending Wilson's theorem (Theorem~\ref{thm:wilson}) to hypergraphs. Moreover, the proof method in \cite{GKLO:ta} is quite different from the one developed in~\cite{keevash:14}, so this gives also a new proof of the Existence conjecture.

\begin{theorem}[\cite{GKLO:ta}] \label{thm:hyperFdec}
	For all $k\in \bN$, $p\in[0,1]$ and any $k$-graph~$F$,
	there exist $c>0$ and $h,n_0\in\bN$ such that the following holds.
	Suppose that $G$ is a $(c,h,p)$-typical $k$-graph on at least $n_0$~vertices. Then $G$ has an $F$-decomposition whenever it is $F$-divisible.
\end{theorem}

Keevash~\cite{keevash:18b} later proved an even more general result. In particular, he obtained a partite version of Theorem~\ref{thm:hyperFdec}. This has the following nice application to \emph{resolvable} decompositions. For simplicity, we only consider graphs. An $F$-decomposition of $G$ is \emph{resolvable} if it can be partitioned into $F$-factors
(where an $F$-factor is a vertex-disjoint union of copies of $F$ covering al vertices
of~$G$).
Recall that Kirkman's famous schoolgirl problem asked for a resolvable 
$K_3$-decomposition of~$K_{15}$.
Now, suppose we seek a resolvable $F$-decomposition of some graph~$G$. 
The number of $F$-factors in such a decomposition should clearly be $t:=\frac{e(G)|F|}{|G|e(F)}$. Define an auxiliary graph $G'$ by adding a set $T$ of $t$ new vertices and all edges between $G$ and~$T$. Moreover, let $F'$ be the graph obtained from $F$ by adding a new vertex $x'$ and all edges to~$V(F)$.
Now, suppose we can find an $F'$-decomposition of $G'$ with the additional property that the copy of $x'$ always lies in~$T$ (\cite{keevash:18b} contains a general framework for this). Then every copy of $F'$ gives a copy of $F$ in~$G$. Moreover, for a fixed vertex in~$T$, the copies of $F'$ containing it induce an $F$-factor. Thus, we have found a resolvable $F$-decomposition of $G$. 
This result, for $F\in \Set{C_3,C_4,C_5}$, was instrumental in the resolution of the Oberwolfach problem (see Section~\ref{sec:OW}).

\subsection{Minimum degree versions}

As Wilson's theorem (Theorem~\ref{thm:wilson}) has been extended to graphs of large enough minimum degree, one can also ask the same for hypergraphs. 
The definition of the decomposition threshold straightforwardly generalises.
For a $k$-graph $F$, define $\delta_F$ as the infimum of all $\delta\in[0,1]$ with the property that any sufficiently large $F$-divisible $k$-graph $G$ with $\delta(G)\ge \delta |G|$ has an $F$-decomposition. 
(Recall that $\delta(G)$ denotes the minimum $(k-1)$-degree of~$G$.)

Explicit bounds for the parameters $c$ and $h$ in Theorem~\ref{thm:hyperFdec} were obtained in~\cite{GKLO:ta}. Since every $k$-graph with sufficiently large minimum degree is typical, this gives explicit upper bounds for~$\delta_F$. However, these bounds are very weak.
This is partly due to a reduction in Theorem~\ref{thm:hyperFdec} that reduces the problem for general $F$ to such $F$ which are `weakly regular'. For weakly regular $F$ (see \cite{GKLO:ta} for the definition), the method in~\cite{GKLO:ta} gives much better results. In particular, we have the following bound on the decomposition threshold of cliques.

\begin{theorem}[\cite{GKLO:ta}] \label{thm:hyper F mindeg}
	For all $k < r$, we have $\delta_{K^k_r} \le 1- \frac{k!}{3\cdot 14^k r^{2k}} $.
\end{theorem}

We remark that even for fractional decompositions, the best known bounds are only slightly better. More precisely (for simplicity we regard the uniformity $k$ as fixed and consider asymptotics in~$r$), it is shown in~\cite{BKLMO:17} that $\delta^*_{K^k_r} \le 1- \Omega_k(r^{-2k+1})$, where $\delta^*_{F}$ is the fractional $F$-decomposition threshold defined in the obvious way.

Note that Theorem~\ref{thm:hyper F mindeg} implies the following result on 
completions of partial Steiner systems: for all $r>k$, there is an $n_0$ so that
whenever $X$ is a `partial' Steiner system, consisting of a set of edge-disjoint  $K^k_r$ on $n$ vertices and $n^*\ge \max \{ n_0,\frac{3\cdot 14^k r^{2k}}{k!}n\}$ satisfies the necessary divisibility conditions, then $X$ can be extended to a $K^k_r$-decomposition of
$K^{k}_{n^*}$. In other words,~$X$ can be extended into a so-called $(n^*,r,k)$-Steiner system.
For the case of Steiner triple systems (i.e.~$r=3$ and $k=2$), Bryant and Horsley~\cite{BH:09} showed that one can take $n^*=2n+1$,
which proved a conjecture of Lindner.
It would be interesting to extend this exact result to other parameter values.

It is not clear what the correct value of $\delta_{K^k_r}$ should be.
As observed in~\cite{GKLO:ta}, a construction from~\cite{KMV:14} can be modified to obtain (for fixed $k<r$) infinitely many $k$-graphs $G$ with $\delta(G) \ge (1- \cO_k( \frac{\log{r}}{r^{k-1}}))|G|$ which are even $K^k_r$-free. Modulo the divisibility of such examples, this seems to suggest that $\delta_{K^k_r}\ge 1- \cO_k( \frac{\log{r}}{r^{k-1}})$. 
In view of the case $k=2$, perhaps the following is true.

\begin{conj}
	$\delta_{K^k_r} = 1- \Theta_k(r^{-k+1})$.
\end{conj}

Moreover, the following generalisation of Corollary~\ref{cor:equal thresholds} might be true.

\begin{conj}
	For all $r>k\ge 2$, $\delta_{K^k_r} = \delta_{K^k_r}^*$.
\end{conj}

Prior to~\cite{GKLO:ta}, the only explicit result for hypergraph decomposition thresholds was due to Yuster~\cite{yuster:00}, who showed that
if $T$ is a linear $k$-uniform hypertree, then every $T$-divisible $k$-graph $G$ on $n$ vertices with minimum vertex degree at least
$(\frac{1}{2^{k-1}}+o(1))\binom{n}{k-1}$ has a $T$-decomposition.
This is asymptotically best possible for nontrivial~$T$. Moreover, the result implies that $\delta_T\le 1/2^{k-1}$.

\section{Euler tours in hypergraphs}\label{sec:euler}

Finding an \emph{Euler tour} in a graph is a problem as old as graph theory itself: Euler's negative resolution of the Seven Bridges of K\"onigsberg problem in 1736 is widely considered the first theorem in graph theory. Euler observed that if a (multi-)graph contains a closed walk which traverses every edge exactly once, then all the vertex degrees are even. Hence, he observed a \emph{necessary divisibility condition} for the existence of such a closed walk, which we now call an \emph{Euler tour}. Is this divisibility condition also sufficient? In general, no, since the graph might be disconnected but still fulfil the divisibility condition. However, if a graph is connected, then it contains an Euler tour if and only if the divisibility condition is satisfied. This fact was already stated by Euler, and its proof is often attributed to Hierholzer and Wiener.

In this section, we consider Euler tours in hypergraphs.
There are several ways of generalising the concept of paths/cycles, and similarly Euler trails/tours, to hypergraphs. We focus here on the `tight' regime. 
Given a $k$-graph $G$, a sequence of vertices $\cW=x_1x_2\dots x_\ell$ is a \emph{(tight self-avoiding) walk in $G$} if $\Set{x_i,x_{i+1},\dots ,x_{i+k-1}}\in E(G)$ for all $i\in [\ell-k+1]$, and no edge of $G$ appears more than once in this way. Similarly, we say that $\cW$ is a \emph{closed walk} if $\Set{x_i,x_{i+1},\dots ,x_{i+k-1}}\in E(G)$ for all $i\in [\ell]$, with indices modulo~$\ell$, and no edge of $G$ appears more than once in this way. We let $E(\cW)$ denote the set of edges appearing in~$\cW$.

\begin{defin}
	An \emph{Euler tour of $G$} is a closed walk $\cW$ in~$G$ with $E(\cW)=E(G)$.
\end{defin}

Clearly, if $G$ is $2$-graph, then this coincides with the usual definition of an Euler tour, and a necessary condition for the existence of such a tour is that every vertex degree is even. Can we formulate an analogous condition for $k$-graphs?
To do so, assume that a $k$-graph $G$ has an Euler tour $\cW$. Fix any vertex $v$. Note that $v$ might appear several times in the sequence $\cW$, however, for every such appearance, it is contained in exactly $k$ edges. Since every edge appears exactly once in $\cW$, we can conclude that the degree $d_G(v)$ of $v$ is divisible by~$k$.
Now, the question is again: is this condition also sufficient for the existence of an Euler tour? Again, the answer is no, as the given hypergraph might be divisible (i.e.~satisfy this degree condition) but consist of disjoint pieces. Moreover, the problem of deciding whether a given $3$-graph has an Euler tour has been shown to be NP-complete~\cite{LNR:17}, thus when $k>2$, there is probably no simple characterisation of $k$-graphs having an Euler tour.
Surprisingly, until recently, it was not even known whether the complete $k$-graph has an Euler tour if it is divisible. Using the language of \emph{universal cycles}, this was formulated as a conjecture by Chung, Diaconis and Graham~\cite{CDG:89,CDG:92} in 1989. More precisely, they conjectured that for every fixed $k\in \bN$ and sufficiently large~$n$, there exists an Euler tour in $K^k_n$ whenever $k$ divides $\binom{n-1}{k-1}$.

Clearly, this is true for $k=2$. 
Numerous partial results have been obtained. In particular, Jackson proved the conjecture for $k=3$~\cite{jackson:93} and for $k\in \Set{4,5}$ (unpublished), and Hurlbert~\cite{hurlbert:94} confirmed the cases $k\in\Set{3,4,6}$ if $n$ and $k$ are coprime (see also~\cite{LTT:81}).
Various approximate versions of the conjecture have been obtained in~\cite{blackburn:12,CHHM:09,DL:16,LTT:81}.

Recently, the conjecture was proven for all $k$ by the authors and Joos~\cite{GJKO:ta}. In fact, the result is more general and applies to quasirandom $k$-graphs in the sense of Definition~\ref{def:typical}. We state a simplified version here which applies for almost complete $k$-graphs.

\begin{theorem}[\cite{GJKO:ta}]\label{thm:Euler}
	For all $k\in \bN$ there exists $\eps>0$ such that any sufficiently large $k$-graph $G$ with $\delta(G)\ge (1-\eps)|G|$ has an Euler tour if $k\mid d_G(v)$ for every $v\in V(G)$.
\end{theorem}

We conjecture that the minimum degree condition can be significantly improved. We discuss this in more detail in Section~\ref{sec:euler problems}.

\subsection{Euler tours: Proof sketch} \label{sec:euler proof}

The proof of Theorem~\ref{thm:Euler} relies on Theorem~\ref{thm:hyperFdec}
in order to complete a suitable partial Euler tour into a `full' one.
More precisely, the proof proceeds as follows. 
We call a walk $\cW$ in a $k$-graph $G$ \emph{spanning} if every ordered $(k-1)$-set of vertices appears consecutively in~$\cW$ at least once. The motivation behind this definition is as follows. Suppose $\cW$ is a closed spanning walk in $G$ and $\cW'$ is some other closed walk which is edge-disjoint from~$\cW$. Then we can `insert' $\cW'$ into $\cW$ as follows: take any $(k-1)$-tuple which appears in $\cW'$. Since $\cW$ is spanning, we know that this ordered tuple appears in $\cW$ too, so we can follow $\cW$ until we reach an appearance of this tuple, then follow $\cW'$ until we reach this tuple again, and then continue with $\cW$. It is easy to see that this yields a new closed spanning walk which uses precisely the edges of $\cW$ and~$\cW'$. Hence, we have the following:
\begin{observe}
	If $G$ can be decomposed into closed walks such that one of them is spanning, then $G$ has an Euler tour.
\end{observe}

This approach breaks the proof into two parts: First, we need to find a spanning walk. 
Note that there are $\Theta(n^{k-1})$ ordered $(k-1)$-sets of vertices, so in order to be spanning, our walk $\cW$ needs to have length $\Omega(n^{k-1})$.  Since $G$ has $\Theta(n^k)$ edges, there is at least enough room. Moreover, if we construct $\cW$ in some random fashion, say using $n^{k-1}\log^{2}n$ edges, then we might hope that every $(k-1)$-set only appears in $\log^{3} n$ edges of~$\cW$. That is, the subgraph formed by the edges of $\cW$ has very small maximum degree, so removing these edges leaves $G$ essentially unchanged.

Now, we come to the second part. We only need to decompose the remainder $G-E(\cW)$ into any number of closed walks. In particular, a decomposition into tight cycles would be sufficient. This we can achieve using the $F$-decomposition result from Section~\ref{sec:Fdecs}, with $F$ being a tight cycle.
Let $C^{k}_\ell$ denote the tight $k$-uniform cycle of length $\ell$, that is, the vertices of $C^{k}_\ell$ are $v_1,\dots,v_\ell$, and the edges are all the $k$-tuples of the form $\Set{v_i,v_{i+1},\dots,v_{i+k-1}}$, with indices modulo~$\ell$.

We want to apply Theorem~\ref{thm:hyperFdec} with $F=C^{k}_{2k}$.
Conveniently, a $k$-graph $G$ is $C^{k}_{2k}$-divisible whenever $2k\mid e(G)$ and $k \mid d_G(v)$ for all $v\in V(G)$, that is, there is no divisibility condition for $i$-sets with $i>1$.\footnote{We have $d_{C^{k}_{2k}}(\Set{v_1,\dots,v_{i-1},v_k})=1$ and hence $d_{C^{k}_{2k}}(i)=1$.} That the vertex degrees in $G-E(\cW)$ are divisible
by $k$ follows automatically from the divisibility of the initial graph $G$ and since we removed a closed walk. Moreover, by removing greedily a few copies of $C^{k}_{2k+1}$, say, we can make the number of (remaining) edges divisible by~$2k$. 
Theorem~\ref{thm:hyperFdec} then does the rest for us.

Let us say a few more words about finding the spanning walk. Essentially, we show that a random self-avoiding walk of length $n^{k-1}\log^{2} n$ has the desired properties with high probability.\footnote{We ignore here that the walk should be closed in the end. This can be easily achieved afterwards in $\cO(1)$ deterministic steps.}
Fix any $k-1$ vertices $X_1,\dots,X_{k-1}$ as a start tuple. Now, in each step, with the current walk being $X_1,\dots,X_{i-1}$, choose a vertex $X_{i}$ uniformly at random from all vertices that form an edge together with the last $k-1$ vertices $X_{i-1},\dots, X_{i-k+1}$ of the current walk and this edge has not been used previously by the walk. If no such vertex exists, then stop.

In order to analyse this random walk, fix vertices $v_1,\dots,v_{k-1}$. We say that the walk \emph{visits} (these vertices) at step $i$ if $X_{i-k+2}=v_1,\dots, X_{i}=v_{k-1}$. 
Now, assume that in some step $i$, the walk visits vertices $v_1',\dots,v'_{k-1}$. We want to ask ourselves: what is the probability that the walk visits $v_1,\dots,v_{k-1}$ at step $i+k$? For simplicity, ignore the condition that the walk ought to be self-avoiding, and also assume that $v_1',\dots,v_{k-1}',v_1,\dots, v_{k-1}$ are distinct.
In each of the following $k$ steps, the walk has clearly at most $n$ choices for the next vertex, so the total number of choices for the walk is at most $n^k$. Crucially, using the minimum degree assumption, we can check that there are $\Omega(n)$ vertices $v^*$ such that $v_1',\dots,v'_{k-1},v^*,v_1,\dots,v_{k-1}$ is a tight walk and thus an admissible choice for the next $k$ steps.
Hence, the (conditional) probability that the walk visits $v_1,\dots,v_{k-1}$ at step $i+k$ is $\Theta(n^{-k+1})$. 
Thus, if we let the walk continue for about $n^{k-1}\log^{2} n$ steps and consider the steps $i$ which are multiples of $k$, then the expected number of visits is roughly $(\log^{2} n)/k$. Moreover, since the stated bound for the probability holds for any outcome of previous such steps, a Chernoff--Hoeffding type inequality applies, and we can infer that the probability of the walk not visiting at all is tiny. A union bound over all ordered $(k-1)$-tuples shows that the walk is spanning with high probability.
A similar argument shows that the walk is unlikely to have maximum degree larger than~$\log^{3} n$. This justifies the above analysis also for the self-avoiding walk, since the number of admissible `link' vertices $v^*$ is still $\Omega(n)$. (Technically, we `stop' the walk as soon as some $(k-1)$-set has too large degree, and then analyse this stopped walk. We refer to~\cite{GJKO:ta} for the remaining details.)

\subsection{Open problems on hypergraph decompositions and Euler tours} \label{sec:euler problems}

The following conjecture would  provide a `genuine' minimum degree version of 
Theorem~\ref{thm:Euler}.

\begin{conj}\label{conj:Euler Dirac}
	For all $k>2$ and $\eps>0$, every sufficiently large $k$-graph $G$ with $\delta(G)\ge (1-\frac{1}{k}+\eps)|G|$ has a tight Euler tour if all vertex degrees are divisible by~$k$.
\end{conj}

It seems possible that the approach for Theorem~\ref{thm:Euler} can be extended 
to attain Conjecture~\ref{conj:Euler Dirac}:
recall that the proof consisted of two steps. First, we found a spanning walk with small maximum degree. For this, we analysed a random self-avoiding walk. The crucial property was that, given any two disjoint $(k-1)$-tuples $(v_1,\ldots,v_{k-1}),(v_{k+1},\ldots,v_{2k-1})$,
there are $\Omega(n)$ vertices $v_k$ such that $v_iv_{i+1}\dots v_{i+k-1}$ is an edge for all $i\in[k]$. This property is already satisfied if $\delta(G)\ge (1-\frac{1}{k} + \eps)|G|$.
In fact, one can even make the argument work if only $\delta(G)\ge (1/2 + \eps)|G|$, by considering longer paths connecting the tuples $(v_1,\ldots,v_{k-1})$ and $(v_{k+1},\ldots,v_{2k-1})$.
Amongst other things, this latter fact motivated the conjecture in~\cite{GJKO:ta} (restated in an earlier version of this survey) that $1/2$ could be the right threshold for Euler tours for any uniformity~$k$.
However, this was apparently too optimistic. Very recently, Piga and Sanhueza-Matamala~\cite{PSM:21} provided a counterexample in the case $k=3$, showing that the threshold needs to be at least~$2/3$. Moreover, they proved that this is the correct threshold (that is, Conjecture~\ref{conj:Euler Dirac} holds for $k=3$).

The bottleneck is the second step, where we decomposed the remaining $k$-graph into tight cycles. In the proof of Theorem~\ref{thm:Euler}, we applied the general $F$-decomposition theorem from Section~\ref{sec:Fdecs} to obtain a decomposition into tight cycles.
The following conjecture\footnote{which replaces the $1/2$ variant from an earlier version that turned out to be false, too, in the case $k=3$, see~\cite{PSM:21}.} would complement the random walk analysis sketched above, and thus imply Conjecture~\ref{conj:Euler Dirac}.
It would also be significant in its own right.
\begin{conj}\label{conj:tight cycle dec}
	Any $k$-graph with $\delta(G) \ge (1-\frac{1}{k}+o(1))|G|$ can be decomposed into tight cycles, provided that all vertex degrees are divisible by~$k$.
\end{conj}

Note that an approximate decomposition is easy to obtain. 
Indeed, since $C^k_{2k}$ is $k$-partite and thus has Tur\'an density~$0$ by a well-known result of Erd\H{o}s~\cite{erdos:64b}, we can iteratively pull out copies of $C^k_{2k}$ until $o(|G|^k)$ edges remain.

\section{Oberwolfach problem}\label{sec:OW}

The Oberwolfach problem, posed by Ringel in 1967, asks for a decomposition of the complete graph $K_{n}$ into edge-disjoint copies of a given $2$-factor. Clearly, this can only be possible if $n$ is odd.

\begin{problem}[Oberwolfach problem, Ringel, 1967]
	Let $n\in \bN$ and let $F$ be a $2$-regular graph on $n$ vertices.
	For which (odd) $n$ and $F$ does $K_{n}$ have an $F$-decomposition?
\end{problem}
The problem is named after the Mathematical Research Institute of Oberwolfach, where Ringel formulated it as follows: assume $n$ conference participants are to be seated around circular tables for $\frac{n-1}{2}$~meals, where the total number of seats is equal to $n$, but the tables may have different sizes. Is it possible to find a seating chart such that every person sits next to any other person exactly once?

Note that when $F$ consists of only one cycle (that is, there is one large table), then we seek a decomposition of $K_n$ into Hamilton cycles, which is possible by Walecki's theorem from 1892. 
In the other extreme, if all tables have only size $3$, then we seek a decomposition of $K_n$ into triangle factors. This was Kirkman's (generalised) schoolgirl problem from 1850, eventually solved by Ray-Chaudhuri and Wilson~\cite{RW:71} and independently by Lu.

Over the years, the Oberwolfach problem and its variants have received enormous attention, with more than 100 research papers produced. 
Most notably, Bryant and Scharaschkin~\cite{BS:09} proved it for infinitely many~$n$.
Traetta~\cite{traetta:13} solved the case when $F$ consists of two cycles only, Alspach, Schellenberg, Stinson and Wagner~\cite{ASSW:89} solved the case when all cycles have equal length, and Hilton and Johnson~\cite{HJ:01} solved the case when all but one cycle have equal length.

An approximate solution to the Oberwolfach problem was obtained by Kim, K\"uhn, Osthus and Tyomkyn~\cite{KKOT:19} and Ferber, Lee and Mousset~\cite{FLM:17}. More precisely, it follows from their (much more general) results that $K_n$ contains $n/2-o(n)$ edge-disjoint copies of any given $2$-factor~$F$.

A related conjecture of Alspach stated that for all odd $n$ the complete graph $K_n$ can be decomposed into any collection of cycles of length at most $n$ whose lengths sum up to~$\binom{n}{2}$.
This was solved by Bryant, Horsley, and Pettersson~\cite{BHP:14}.

Very recently, the Oberwolfach problem was solved by the authors together with Joos and Kim~\cite{GJKKO:18}. More precisely, they showed that for all odd $n\ge n_0$, there is a solution for any given $2$-factor~$F$. The remaining cases could (in theory) be decided by exhaustive search, but this is not practically possible as $n_0$ is rather large. 
It would be very interesting to complete the picture. Perhaps there are not many exceptions (currently, there are four known exceptions).

We state the result in the following slightly more general way, where $K_n$ can be replaced by an almost-complete graph. Recall that this means one can obtain a solution to the Oberwolfach problem even if the first $o(n)$ copies of $F$ are chosen greedily.

\begin{theorem}[\cite{GJKKO:18}]\label{thm:OW}
	There exists $\eps>0$ such that for all sufficiently large $n$, the following holds: Let $F$ be any $2$-regular graph on $n$ vertices, and let $G$ be a $d$-regular graph on $n$ vertices for some even $d\in \bN$.
	If $d \ge (1-\eps)n$, then $G$ has an $F$-decomposition.
\end{theorem}
As mentioned earlier, the proof relies on Corollary~\ref{cor:OW mindeg approx}
(to obtain a suitable approximate decomposition)
and the results on resolvable cycle decompositions in~\cite{keevash:18b}
(as part of an absorbing approach).

In the spirit of this survey, the obvious question is of course: can the minimum degree assumption in Theorem~\ref{thm:OW} be weakened? We discuss this further in Section~\ref{sec:OW problems}.

An immediate consequence of Theorem~\ref{thm:OW} is that if $n$ is even, then $K_n$ can be decomposed into one perfect matching and otherwise copies of~$F$. For this variant of the Oberwolfach problem as well, many partial results were previously obtained (see e.g.~\cite{BD:11,HS:91,HKR:79}).

Another natural extension is the following `generalised Oberwolfach problem'. Suppose $F_1,\dots,F_{(n-1)/2}$ are (possibly distinct) $2$-factors on $n$ vertices. Is it possible to decompose $K_n$ into $F_1,\dots,F_{(n-1)/2}$? 
In the special case where the list contains only two distinct $2$-factors, this is known as the Hamilton--Waterloo problem, which was also solved in~\cite{GJKKO:18} (for sufficiently large~$n$). In fact, Theorem~\ref{thm:OW} holds in this general setting provided that some $2$-factor appears linearly many times in the list. 

Improving on this, Keevash and Staden~\cite{KS:20a} recently solved the generalised Oberwolfach problem. Their result applies in the setting of dense typical graphs, and they also prove an appropriate version of this for directed graphs. 

\begin{theorem}[\cite{KS:20a}]
	For every $p>0$ there exist $c>0$ and $h\in \bN$ such that the following holds. Any sufficiently large $(c,h,p)$-typical graph $G$ which is $d$-regular for some even $d\in \bN$ can be decomposed into any $d/2$ given $2$-factors.
\end{theorem}

\subsection{Open problems related to the Oberwolfach problem}\label{sec:OW problems}

\subsubsection{Minimum degree thresholds}

As we have already seen, the problem of decomposing a graph $G$ into a given $2$-factor not only makes sense if $G$ is complete. Of course, we should assume that $G$ is regular with even degrees.

\begin{conj}[\cite{GJKKO:18}]
	For all $\eps>0$, the following holds for sufficiently large~$n$.
	Let $F$ be any $2$-regular graph on $n$ vertices, and let $G$ be a $d$-regular graph on $n$ vertices for some even $d\in \bN$.
	If $d \ge (3/4+\eps)n$, then $G$ has an $F$-decomposition.
\end{conj}

If $F$ is a triangle-factor, then the `threshold' $3/4$ would be optimal. On the other hand, if $F$ is a Hamilton cycle, then it can be lowered to~$1/2$ (\cite{CBKLOT:16}). It would be interesting to `interpolate' between these extremal cases.
More specifically, it could be true that if all cycles in $F$ have even length, then $2/3$ is sufficient, and if $C_4$ is excluded in addition, then the threshold is~$1/2$. Similarly, if the girth of $F$ 
is sufficiently large compared to $1/\eps$, then~$1/2$ should also be sufficient.
Note that Corollaries~\ref{cor:OW mindeg approx} and~\ref{cor:OW mindeg approx girth} give some partial approximate results in this direction.

\subsubsection{Hypergraphs}

It seems natural to ask for an analogue of the Oberwolfach problem for hypergraphs. Bailey and Stevens~\cite{BS:10} conjectured that $K^k_n$ has a decomposition into tight Hamilton cycles if and only if $k$ divides $\binom{n-1}{k-1}$. This is still open and would generalise Walecki's theorem to hypergraphs.
Clearly, the condition $k \mid \binom{n-1}{k-1}$ is necessary since every Hamilton cycle contains $k$ edges at any fixed vertex. Moreover, it also implies that the total number of edges of $K^k_n$ is divisible by $n$, the number of edges in one Hamilton cycle. We conjecture that the same divisibility condition guarantees a decomposition into any tight cycle factor.

\begin{conj}
	For fixed $k$, the following holds for sufficiently large~$n$.
	Let $F$ be the vertex-disjoint union of tight $k$-uniform cycles, each of length at least $2k-1$, with $n$ vertices in total.
	Then $K^k_n$ has an $F$-decomposition if $k$ divides $\binom{n-1}{k-1}$.
\end{conj}

To the best of our knowledge, this has not been explicitly asked before. We note that the somewhat generous assumption that each cycle has length at least $2k-1$ ensures that there are no divisibility obstructions for $i$-sets with $i>1$.
We also note that the very general result of Ehard and Joos~\cite{EJ:20b} concerning approximate decompositions of quasirandom hypergraphs into bounded degree subgraphs yields an approximate solution to the above conjecture, in that $K^k_n$ contains any collection of $(1-o(1))\binom{n-1}{k-1}/k$ edge-disjoint tight cycle factors.

Another problem which is related to the conjecture of 
Bailey and Stevens was made by Baranyai as well as Katona. First, recall Baranyai's theorem~\cite{baranyai:75} stating that
$K^k_n$ has a $1$-factorization whenever $k$ divides~$n$. 
As in the case of graphs, a \emph{$1$-factor} (or \emph{perfect matching}) is a set of disjoint edges covering all the vertices,
and a \emph{$1$-factorization} is a set of edge-disjoint $1$-factors covering all the edges.
Baranyai~\cite{baranyai:79} and Katona conjectured an extension of 
Baranyai's theorem to the case when the divisibility condition $k \mid n$ is not satisfied:
instead of decomposing into $1$-factors, the aim is to decompose into \emph{wreaths}.
Here, given a cyclically ordered set of $n$ vertices, a wreath is obtained by greedily choosing hyperedges as follows: the first hyperedge consists
of $k$ consecutive vertices and in each step the next hyperedge consists of the 
$k$ vertices which come directly after the vertices in the previous hyperedge. This process stops as soon as one obtains a regular hypergraph.
So for instance, if $k=4$ and $n=6$, then $\{1234, 5612, 3456 \}$ is a wreath.
If $k \mid n$, then a wreath is a perfect matching.
If $n$ and $k$ are co-prime, then a wreath is a tight Hamilton cycle.
The \emph{wreath decomposition conjecture} postulates that $K^k_n$ can always be decomposed into wreaths.

While this problem is still open for the complete hypergraph, we propose the following minimum degree version of Baranyai's theorem.

\begin{conj}
	For fixed $k$ and $\eps>0$, the following holds for all sufficiently large~$n$.
	An $n$-vertex $k$-graph $G$ with $\delta(G)\ge (1/2+\eps)n$ can be decomposed into perfect matchings if and only if $k\mid n$ and $G$ is vertex-regular.
\end{conj}

\subsubsection{Decompositions into \texorpdfstring{$r$}{r}-factors}

Finally, we restate the following conjecture formulated in~\cite{GJKKO:18}, which 
can be viewed as a far-reaching generalisation of the Oberwolfach problem
from $2$-regular to regular graphs of arbitrary degrees.

\begin{conj}[\cite{GJKKO:18}]
	For all $\Delta \in \mathbb{N}$, there exists an $n_0 \in \mathbb{N}$ so that the following holds for all $n \ge n_0$.
	Let $F_1,\ldots,F_t$ be $n$-vertex graphs such that $F_i$ is $r_i$-regular 
	for some $r_i \le \Delta$ and $\sum_{i\in [t]}r_i={n-1}$.
	Then there is a decomposition of $K_n$ into $F_1,\ldots,F_t$.
\end{conj}

This conjecture is clearly extremely challenging.
So it would be interesting to prove it for restricted families, such as graphs which 
are separable or have high girth.
An approximate version of the above conjecture was proved
by Kim, K\"uhn, Osthus and Tyomkyn~\cite{KKOT:19}.

\section{Related decomposition problems}\label{sec:related}

In this final section, we briefly mention some further decomposition problems. We also remark that, for all the questions we discussed, it is interesting to ask for algorithmic variants (can a decomposition of a dense hypergraph be found in polynomial time?), counting problems (how many different decompositions of a graph exist?) and many other directions which we did not cover here.

\subsection{Weighted decompositions into triangles and edges}
Recall Kirkman's theorem that $K_n$ has a triangle decomposition whenever it is divisible.
We might ask what happens if we ignore divisibility? Can we decompose into triangles and a few edges?
For a decomposition of an $n$-vertex graph $G$ into $e$ edges and $t$ triangles, we define $2e+3t$ as the weight of the decomposition, and the aim is to find a decomposition of minimum weight, denoted by $\pi_3(G)$.\footnote{More generally, suppose a fixed set $\cH$ of graphs and a weight function $w$ on $\cH$ are given. For a graph $G$ which is decomposed into $H_1,\dots,H_s\in \cH$, define the weight of this decomposition as $\sum_{i=1}^s w(H_i)$. One can then ask for the minimum weight of such a decomposition of $G$. 
	The case when $\cH$ is the set of all cliques and $w(K_r)=r$ has received much attention.}
Clearly, we always have $\pi_3(G)\ge e(G)$. In particular, $$\pi_3(K_n)\ge \binom{n}{2} = (1+o(1))n^2/2.$$
Similarly, if $G$ is triangle-free, then $\pi_3(G)= 2e(G)$. In particular, $$\pi_3(K_{\lceil n/2\rceil,\lfloor n/2\rfloor})= 2\cdot \lceil n/2\rceil \cdot \lfloor n/2\rfloor = (1+o(1))n^2/2 .$$

Define $\pi_3(n)$ as the maximum of $\pi_3(G)$ over all $n$-vertex graphs~$G$.
The problem of determining $\pi_3(n)$ was first considered by Gy\"ori and Tuza~\cite{GT:87}.
Kr\'{a}l', Lidick\'{y}, Martins and Pehova~\cite{KLMP:19} resolved this problem
asymptotically, by showing that $$\pi_3(n)= (1/2+o(1))n^2.$$
Blumenthal, Lidick\'{y}, Pehova, Pfender, Pikhurko and Volec~\cite{BLPPPV:19}
were able to strengthen this to an exact bound. It turns out that $K_n$ and $K_{\lceil n/2\rceil,\lfloor n/2\rfloor}$ are the only extremal examples.
A crucial tool in the proof of these results was the triangle case of Theorem~\ref{thm:dec main upper}, proved in~\cite{BKLO:16}.

An immediate consequence of the above results is that every $n$-vertex graph with
$n^2/4 +k$ edges contains $2k/3-o(n^2)$ edge-disjoint triangles.
A problem of Tuza~\cite{tuza:01} would generalise the latter bound to arbitrary cliques.
\begin{problem}[\cite{tuza:01}]
	Does every $n$-vertex graph with $\frac{r-2}{2r-2}n^2+k$ edges contain $\frac{2}{r}k-o(n^2)$ edge-disjoint copies of $K_r$?
\end{problem}

A minimum degree version of this problem was considered by Yuster~\cite{yuster:14}: what is the largest number of edge-disjoint copies of $K_r$ one can find in a graph of given minimum degree? Again, this question is still open.

\subsection{Packing and covering number}

For graphs $F$ and $G$, the \emph{$F$-packing number} of $G$, denoted $P(F,G)$, is the maximum number of edge-disjoint copies of $F$ in~$G$. 
The `dual' notion is the \emph{$F$-covering number}, denoted $C(F,G)$, which is the minimum number of copies of $F$ in $G$ that cover every edge of~$G$. Clearly, $G$ has an $F$-decomposition if and only if $P(F,G)=C(F,G)=e(G)/e(F)$.

Assume $F$ is a fixed graph and $G$ is a sufficiently large and dense graph. Note that a collection of edge-disjoint copies of $F$ in $G$ forms an $F$-decomposition of some subgraph of~$G$. Hence, a natural way for proving lower bounds on $P(F,G)$ is to delete as few edges as possible from $G$ to obtain an $F$-divisible graph, and then to apply a decomposition result such as Gustavsson's theorem (see Section~\ref{sec:thresholds}). Using this approach, Caro and Yuster~\cite{CY:97,CY:98} determined $P(F,K_n)$ and $C(F,K_n)$ exactly for all sufficiently large~$n$.
Moreover, Alon, Caro and Yuster~\cite{ACY:98} proved that when $G$ is large and very dense ($\delta(G)\ge (1-o(1))|G|$), then $P(F,G)$ and $C(F,G)$ can be computed in polynomial time.
The minimum degree threshold in this result can probably be significantly improved, perhaps a natural guess is that $\delta(G)\ge (\delta_F+o(1))|G|$ suffices.

\subsection{Decomposing highly connected graphs into trees}

Essentially all decomposition results we discussed in this survey apply only for dense graphs with linear minimum degree.
It would be very interesting to investigate different conditions which ensure that a given $F$-divisible graph $G$ has an $F$-decomposition. 

One such example is a beautiful conjecture of Bar\'{a}t and Thomassen~\cite{BT:06} on decompositions into a fixed tree~$T$. Recall from Section~\ref{sec:bip dec} that the decomposition threshold of $T$ is $1/2$. Moreover, since $\gcd(T)=1$, the only necessary divisibility condition for $G$ to have a $T$-decomposition is $e(T)\mid e(G)$.
The reason why the minimum degree threshold cannot be lowered is that $G$ could consist of two equal-sized vertex-disjoint cliques (with a few edges removed), such that the total number of edges is divisible by $e(T)$, but the number of edges in each clique is not.
However, this example is not very robust. Just adding a constant number of edges across would allow us to find a $T$-decomposition. 
In particular, if $G$ is highly connected, then it seems hard to construct any such example. Bar\'{a}t and Thomassen~\cite{BT:06} conjectured that in fact this is impossible. This was proved recently by Bensmail, Harutyunyan, Le, Merker and Thomass\'{e}~\cite{BHLMT:17} via probabilistic methods, but also 
involving tools based on nowhere-zero flows~\cite{merker:17,thomassen:13}. 

\begin{theorem}[\cite{BHLMT:17}]\label{thm:BT}
	For any tree $T$, there exists a constant $k_T$ such that any graph $G$ which is $k_T$-edge-connected and satisfies $e(T)\mid e(G)$ has a $T$-decomposition.
\end{theorem}

The value of $k_T$ needed for their proof is quite large, and it would be interesting to improve it.

\subsection{Tree packings}\label{sec:trees}

We now discuss some further results on decompositions into trees. Since the term `tree decomposition' is already reserved for another graph-theoretical concept, this problem is usually referred to as `tree packing'.
The main open problem in the area is the so-called `tree packing conjecture' due to 
Gy\'{a}rf\'{a}s and Lehel~\cite{GL:78}.

\begin{conj}[\cite{GL:78}]\label{conj:tree packing}
	For every $n$, the complete graph $K_n$ can be decomposed into any sequence of trees $T_1,\dots,T_n$ where $|T_i|=i$.
\end{conj}

Joos, Kim, K\"uhn and Osthus~\cite{JKKO:19} proved this for bounded degree trees.
Slightly earlier, Allen, B\"ottcher, Hladky and Piguet~\cite{ABHP:19} proved an approximate version
for trees whose maximum degree is allowed to be as large as $o(n/\log n)$.
Very recently, Allen, B\"ottcher, Clemens and Taraz~\cite{ABCT:19} showed that the tree packing conjecture holds for almost all sequences of trees.
Each of these results applies in a considerably more general setting than stated here,
and there are many more results which we do not mention here.

Another famous question on tree packings was formulated by Ringel in 1963, who asked whether $K_{2n+1}$ can be decomposed into any tree with~$n$ edges. Very recently, Montgomery, Pokrovskiy and Sudakov~\cite{MPS:20} and
Keevash and Staden~\cite{KS:20b}
solved Ringel's conjecture for large enough~$n$.

\begin{theorem}[\cite{KS:20b,MPS:20}]
	For sufficiently large~$n$, $K_{2n+1}$ can be decomposed into any tree with~$n$ edges.
\end{theorem}

The proof in~\cite{MPS:20} is based on finding a single rainbow copy of the desired tree~$T$
in a suitably edge coloured $K_{2n+1}$.
The approach in~\cite{KS:20b} builds on results in~\cite{keevash:18b}. One crucial ingredient in both papers is to consider three cases according to the structure of the given tree, which was developed in~\cite{MPS:ta} to prove an approximate version of Ringel's conjecture.

The following conjecture of Graham and H\"{a}ggkvist~\cite{haggkvist:89}  
generalises Ringel's conjecture to arbitrary regular graphs.

\begin{conj}[\cite{haggkvist:89}]\label{conj:GH}
	For any tree $T$, any $2e(T)$-regular graph $G$ has a $T$-decom-position.
\end{conj}

The main result in~\cite{KS:20b} implies that this is true if the host graph $G$ is in addition dense and quasirandom.
Moreover, Corollary~\ref{cor:treebandwidth} gives an approximate version
if $|T| \ge (1/4+o(1))|G|$ and $T$ has bounded maximum degree.
It would be interesting to settle this `dense' case exactly.

\subsection{Sparse decompositions of dense graphs: Erd\texorpdfstring{\H{o}}{\"o}s meets Nash-Williams}

Recall Kirkman's theorem that every $K_3$-divisible complete graph has a $K_3$-decom-position. 
Much of the content of this survey has been inspired by Conjecture~\ref{conj:NW}, which would be a far-reaching generalisation of Kirkman's theorem, and which was posed by Nash-Williams in 1970. Around the same time, Erd\H{o}s proposed another beautiful extension of Kirkman's theorem.
Define the \emph{girth} of a set $\cT$ of triangles to be the smallest $g\ge 4$ for which there is a set of $g$ vertices spanning at least $g-2$ triangles from the set $\cT$. Note that any $K_3$-decomposition has girth at least~$6$.
Erd\H{o}s~\cite{erdos:73} conjectured that there are Steiner triple systems (i.e.~$K_3$-decompositions of $K_n$) of arbitrarily large girth. (Decompositions with high girth are also called `locally sparse' since any set of $4\le j< g$ vertices contains at most $j-3$ triangles.)

\begin{conj}[\cite{erdos:73}]\label{conj:erdos}
	For every fixed $g$, any sufficiently large $K_3$-divisible $K_n$ has a $K_3$-decomposition with girth at least~$g$.
\end{conj}

This conjecture has been proved exactly only for the first non-trivial case, namely $g=7$, in a series of papers~\cite{brouwer:77,GGW:00,GMP:90,LCGG:00}. Recently, it was solved approximately for all fixed~$g$~\cite{BW:19,GKLO:20}.
A generalisation of Conjecture~\ref{conj:erdos} to Steiner systems with arbitrary parameters was formulated in~\cite{GKLO:20}.

We are tempted to propose the following combination of the conjectures of Erd\H{o}s and Nash-Williams.

\begin{conj}
	For every fixed $g$, any sufficiently large $K_3$-divisible graph $G$ with $\delta(G)\ge 3|G|/4$ has a $K_3$-decomposition with girth at least~$g$.
\end{conj}

Of course, given that both conjectures themselves are still open, this seems very challenging. In view of this,
it would even be interesting to obtain approximate decompositions of large girth in a sufficiently dense graph, for instance, to show that any sufficiently large graph $G$ with $\delta(G)\ge 0.9|G|$, say, has an approximate $K_3$-decomposition with arbitrarily high girth.


\providecommand{\bysame}{\leavevmode\hbox to3em{\hrulefill}\thinspace}
\providecommand{\MR}{\relax\ifhmode\unskip\space\fi MR }
\providecommand{\MRhref}[2]{%
	\href{http://www.ams.org/mathscinet-getitem?mr=#1}{#2}
}
\providecommand{\href}[2]{#2}

\end{document}